\documentclass[12pt]{amsart}
\usepackage{times}
\usepackage{amsfonts}
\usepackage{amsthm}
\usepackage{amsmath}
\advance \topmargin by -\headheight
\advance \topmargin by -\headsep
\evensidemargin \oddsidemargin
\newtheorem{theorem}{Theorem}[section]
\newtheorem{lemma}[theorem]{Lemma}
\newtheorem{corollary}[theorem]{Corollary}

\begin{document}

\title{Dimension and Entropy Computations for $L(F_r)$ }

\author{Kenley Jung}

\address{Department of Mathematics, University of California,
Berkeley, CA 94720-3840,USA}

\email{factor@math.berkeley.edu}
\thanks{Research supported by the NSF Graduate Fellowship Program}

\begin{abstract}  We show that certain generating sets of Dykema and 
Radulescu for $L(F_r)$ have free Hausdorff dimension $r$ and nondegenerate 
free Hausdorff $r$-entropy.  
\end{abstract}

\maketitle

      In [6] Voiculescu showed that certain compressions of free group
factors are again free group factors. Radulescu generalized this in [4]
where he introduced the interpolated free group factors denoted by
$L(F_r), r \geq 1.$ These $II_1$-factors coincide with the free group
factors for integral values of $r.$ Dykema ([1]) independently discovered
such factors.  In addition to extending the compression formula of [6],
[1] and [4] showed that the interpolated free group factors are either all
mutually isomorphic or they are all mutually nonisomorphic.

\vspace{.1in}

Voiculescu also introduced in [7] the notion of a microstate as part of
his working theory of free probability.  Microstates allow one to make
sense of Lebesgue measure and Minkowski dimension of $n$-tuples in a
tracial von Neumann algebra.  Using microstates [3] took a fractal
geometric approach by introducing free Hausdorff $r$-entropy, denoted by
$\mathbb H^r.$ It is a kind of asymptotic Hausdorff measurement on the
microstate spaces.  For integral values, free Hausdorff entropy is a
normalization of free entropy:  $\mathbb H^n(z_1,\ldots,z_n) =
\chi(z_1,\ldots,z_n) + \frac{n}{2} \log (\frac{2n}{\pi e}).$ Thus one can
view $\mathbb H^r$ as a continuous extension of $\chi.$ In view of the
interpolated free group factors, a natural question is whether $\mathbb
H^r$ has the same relationship to $L(F_r)$ that $\chi$ has to $L(F_n).$ We
simply mean the following.  It was shown early on ([7]) that there exist
$n$ self-adjoint generators $s_1, \ldots, s_n$ for $L(F_n)$ for which $-
\infty < \chi(s_1,\ldots,s_n) < \infty.$ We want to know whether there are
generators $z_1,\ldots,z_n$ for the interpolated free group factor
$L(F_r)$ satisfying

\[ -\infty < \mathbb H^r(z_1,\ldots,z_n) < \infty.\]

\vspace{.1in}

This remark verifies that certain kinds of generators of Dykema and
Radulescu satisfy the above inequality.  For a larger class of generators
we are able to show that their free Hausdorff dimension (denoted by
$\mathbb H$) is $r.$ This is weaker for $- \infty < \mathbb
H^r(z_1,\ldots,z_n) < \infty \Rightarrow \mathbb H(z_1,\ldots z_n) = r$ in
the same way that $-\infty < \chi(z_1,\ldots,z_n) < \infty \Rightarrow
\delta_0(z_1,\ldots,z_n)=n.$

\vspace{.1in}

More specifically we show the dimension equation for finitely many of
Radulescu's generators along with their natural projections.  We also show
this for finitely many generators (with some natural projections) which
correspond to Dykema's definition of $L(F_r)$ but we must restrict
ourselves to commuting projections due to computational limitations.  If
we assume that the join of the natural projections in either of the two
situations is strictly dominated by the identity operator, then we also
get finite free Hausdorff entropy for the generators.  The estimates are
the usual kind in free dimension approximations.  The generators are
obtained from a combination of three types of algebraic properties:
freeness, commutativity, and finite-dimensionality.  Free Hausdorff
entropy and dimension, though not as user friendly as $\chi$ and
$\delta_0,$ can still cope with such relations and provide the right
bounds.

\vspace{.1in}

Section 1 reviews the definitions of $L(F_r)$ given by Dykema and 
Radulescu and makes precise the kinds of sets we will estimate.  
Sections 2 and 3 demonstrate the upper and lower bounds, respectively.

\section{Review and Notation}

In this brief section we recall the results of Radulescu and Dykema
concerning the interpolated free group factors and make assumptions to be
held for the remainder of this paper.

\vspace{.07in}

$M_k(\mathbb C)$ denotes the set of $k \times k$ matrices over $\mathbb
C,$ $M^{sa}_k(\mathbb C)$ denotes the set of $k \times k$ self-adjoint
matrices over $\mathbb C$ and $tr_k$ is the tracial state on $M_k(\mathbb
C).$ $U_k$ is the group of $k \times k$ unitaries. For $p \in \mathbb N$
$(M_k(\mathbb C))^p$ and$(M^{sa}_k(\mathbb C))^p$ are the spaces of
$p$-tuples of elements of $M_k(\mathbb C)$ or $M^{sa}_k(\mathbb C),$
respectively.  For any subset $X$ of $(M_k(\mathbb C))^p$ or
$(M_k^{sa}(\mathbb C))^p$ and $u \in U_k$ write $uXu^*$ for the subset
obtained by conjugating each entry of an element of $X$ by $u.$ For any $p
\in \mathbb N$ $|\cdot|_2$ is the norm on $(M_k(\mathbb C))^p$ given by
$|(x_1, \ldots, x_p)|_2 = \left ( \sum_{i=1}^p tr_k (x_i^*x_i) \right
)^{\frac{1}{2}}$ and $\|(x_1,\ldots, x_p)\|_2 = \sqrt{k}
|(x_1,\ldots,x_p)|_2.$ $|\cdot|_{\infty}$ will denote the operator norm.  
We maintain the notation for $\Gamma(:), P_{\epsilon}, K_{\epsilon},
\delta_0, \mathbb H^r, \mathbb P^r$ and all other quantities introduced in
[3], [7], and [8] until otherwise stated (see the remark on $M_k(\mathbb
C)$ and $M^{sa}_k(\mathbb C)$ microstates below).

\vspace{.1in}

$\{s\}\bigcup <s_i>_{i=1}^{n},$ $n \in \mathbb N \bigcup \{\infty\}$ is a
family of free semicircular elements in a von Neumann algebra $M$ with
normal, tracial state $\varphi$ and identity $I.$ Suppose $<e_i,
f_i>_{i=1}^{n},$ is a family of nonzero orthogonal projections in
$\{s\}^{\prime \prime}$ such that for any $i$ either $e_i =f_i$ or $e_i
f_i = 0.$ Set $r = 1 + \sum_{i=1}^{n} m_i \alpha_i \beta_i$ where
$\alpha_i = \varphi(e_i), \beta_i = \varphi(f_i),$ $m_i =2$ if $e_i$ and
$f_i$ are orthogonal, and $m_i = 1$ if $e_i = f_i.$ If the von Neumann
subalgebra of $M$ generated by $s$ and $<e_is_if_i>_{i=1}^{n}$ is a
factor, then Radulescu called this the interpolated free group factor,
$L(F_r).$ He showed that this definition made sense, i.e., it depends not
upon the choices of the projections but only upon the value $r$ defined in
terms of the traces of the projections.

\vspace{.07in}

Dykema approached $L(F_r)$ in a slightly different manner.  Suppose
$\mathcal R$ is a copy of the hyperfinite $II_1$-factor free with respect
to $<s_i>_{i=1}^n.$ Given any family of projections $<e_i>_{i=1}^n$ in
$\mathcal R$ Dykema defined $L(F_r)$ for $r >1$ to be the von Neumann
algebra generated by $\mathcal R$ and $<e_i s_i e_i>_{i=1}^n$ where $r = 1
+ \sum_{i=1}^n \varphi(e_i)^2.$ It is a consequence of his work in [1]
that this definition also makes sense and that any such sets along with
$\mathcal R$ generate a factor.

\vspace{.07in}

	We will show some of these generators along with their associated
projections have free Hausdorff dimension $r$ and in the case where the
join of the $e_i$ and $f_i$ is strictly dominated by the identity, we will
show that these generating sets will have finite free Hausdorff
$r$-entropy.  For this we make a few remarks.  We only consider
free dimensions and entropies for finite sets of elements.  For the
remainder of this paper we assume $n \in \mathbb N$ and $<e_i,
f_i>_{i=1}^n$ are projections in $\{s\}^{\prime \prime}$ such that for any
$i$ either $e_i = f_i$ or $e_if_i =0.$ $r = 1 + \sum_{i=1}^n m_i \alpha_i
\beta_i$ as discussed above in the paragraph on Radulescu's generators.  
$\mathcal R$ is a copy of the hyperfinite $II_1$-factor free with respect 
to $<s_i>_{i=1}^N$ and $s \in \mathcal R.$  Also, if for some $i$ $e_i$ 
and $f_i$ are orthogonal, then $\{s,e_1s_1f_1, \ldots, e_ns_nf_n\}$ consists of non self-adjoint elements.  
Thus, we must appropriately modify our definition of $\mathbb H$ and
$\mathbb H^r.$ This is easily done.  We simply take Definition 3.2 of [3]
and replace the $M^{sa}_k(\mathbb C)$ microstates of a self-adjoint
$n$-tuple with the $M_k(\mathbb C)$ microstates for an (not necessarily
self-adjoint)  $n$-tuple (here the $M_k(\mathbb C)$ microstates
approximate the $n$-tuple in $*$-moments).  For the remainder of the paper
microstate spaces and quantities will be taken with respect to
$M_k(\mathbb C)$ approximants and if the quantities are taken with respect
to self-adjoint microstates, then we distinguish them from the
non-self-adjoint quantities with $^{sa}$ (e.g.  $\Gamma_R^{sa}, \chi^{sa},
\mathbb H^{r, sa},$ $\mathbb P^{\alpha, sa},$ etc.).

\vspace{.07in}

Given finite ordered sets $X_1= \{x_{11}, \ldots,
x_{1n_1}\},\ldots, X_p= \{x_{p1},\ldots,x_{pn_p}\}$ of 
elements in $M$ we write $\Gamma_R(X_1,\ldots,X_p;m,k,\gamma)$ for \[
\Gamma_R( x_{11}, \ldots, x_{1n_1},\ldots,
x_{p1},\ldots,x_{pn_p};m,k,\gamma) \] Similarly we will abbreviate all
associated entropies and dimensions.  Set $E = \{e_1,\ldots,e_n\}, F =
\{f_1,\ldots,f_n\},$ and $G = \{e_1 s_1 f_1,\ldots, e_n s_n f_n\}.$ $B$
will denote any $l$-tuple ($l \in \mathbb N$) of strictly contractive
self-adjoint elements in $\mathcal R.$ We do not exclude the situation 
where $B =\emptyset.$  We assume that for some $i$ both $e_i$ and $f_i$ 
are nonzero.

\vspace{.07in}  

Finally, our goal is to show that $\mathbb H(s,B,E,F,G) =r.$ Moreoever,
we show that if $\vee_{i=1}^n (e_i \vee f_i) < I,$ then $- \infty < 
\mathbb
H^r(s,B,E,F,G) < \infty.$ If $B = \emptyset,$ then we have generators for
$L(F_r)$ of Radulescu's type and if $B$ is chosen so that $s$ and $B$
generate $\mathcal R$ and $e_i =f_i$ for each $i,$ then we have generators
that fall into Dykema's picture of $L(F_r).$ We point out that Dykema's
generators are self-adjoint and can be considered with self-adjoint
quantities.  The arguments of Section 2 and 3 work equally well in the
self-adjoint context to show the corresponding statements for the 
self-adjoint entropies and dimensions.

\section{Upper Bound}

\begin{lemma} $\mathbb H^r(s,B,E,F,G) \leq \mathbb P^r(s,B,E,F,G) < 
\infty$ \end{lemma}

\begin{proof} By [3] it suffices to show that $\mathbb P^r(s,B,E,F,G) <
\infty.$ By [5] there exist $C, \epsilon_0 >0$ such that for each $0 <
\epsilon < \epsilon_0$ there exists an $\epsilon$-cover $<u^{(k)}_j>_{j 
\in J_k}$ of $U_k$ taken with respect to $|\cdot|_{\infty}$ with $\# J_k <
\left (\frac{C}{\epsilon} \right)^{k^2}.$ Suppose $R>0$ and $0 < \epsilon
< \epsilon_0$.  $\{s\} \bigcup B \bigcup E \bigcup F$ generate a
hyperfinite von Neumann algebra and thus by Lemma 4.2 of [2] there exist
$m \in \mathbb N$ and $\gamma >0$ such that if $(x,T,P,Q), (x^{\prime},
T^{\prime}, P^{\prime},Q^{\prime}) \in \Gamma_R(s,B, E,F;m,k,\gamma),$
then there exists a $u \in U_k$ such that

\[ |u(x,T,P,Q,)u^* - (x^{\prime}, T^{\prime},P^{\prime}, Q^{\prime})|_2 
<\frac {\epsilon}{(R+1)^2}.\]
 
\noindent We point out that by choosing $m$ and $\gamma$ appropriately the
inequality holds in our non self-adjoint context above, even though we are
using $M_k(\mathbb C)$ microstates.  This is because all the operators in
question are self-adjoint and thus the problem reduces to the self-adjoint
situation.

\vspace{.1in}

For sufficiently large $k$ fix $(x^{(k)}, B^{(k)}, E^{(k)}, F^{(k)}) \in
\Gamma^{sa}_1(s,B,E,F;m,k,\gamma).$ We can arrange it so that writing
$E^{(k)} = \{e_1^{(k)}, \ldots, e_n^{(k)} \}$ and $F^{(k)} = \{f_1^{(k)},
\ldots, f_n^{(k)} \}$ for each $i$ $e_i^{(k)}$ and $f_i^{(k)}$ are
projections and $tr_k(e_i^{(k)}) \leq \alpha_i$ and $tr_k(f_i^{(k)}) \leq
\beta_i.$

\vspace{.1in}

If $e_if_i = 0$ then $\Omega_i^{(k)}$ denotes the ball of $| \cdot |_2$
radius $2$ in $M_k(\mathbb C)$ centered at the origin and if $e_i =f_i$
then $\Omega_i^{(k)}$ denotes the ball of $| \cdot |_2$ radius $2$ in
$M^{sa}_k(\mathbb C)$ centered at the origin.  For each $i$ and $k$
$e_i^{(k)} \Omega_i^{(k)} f_i^{(k)}$ is isometric (when endowed with the
$|\cdot|_2$ metric) to a ball of radius $2$ in Euclidan space of dimension
no greater than $m_i \alpha_i \beta_i.$ For each $i$ and $k$ we can find
an $\epsilon$-cover (with respect to $| \cdot |_2$ ) for $e_i^{(k)}
\Omega_i^{(k)} f_i^{(k)}$ with cardinality no greater than

\[ \left ( \frac{3}{\epsilon} \right)^{m_i \alpha_i \beta_i k^2}.
\]

\noindent Hence we can find an $3\epsilon \sqrt{n}$-net 
$<G_h^{(k)}>_{h \in H_k}$ with respect to 
$| \cdot|_2$ for $e_1^{(k)} \Omega_1^{(k)} f_1^{(k)} \times \cdots \times
e_n^{(k)} \Omega_n^{(k)} f_n^{(k)}$ satisfying

\[ \# H_k \leq \left( \frac{3}{\epsilon} \right)^{(r-1)k^2}.\]

\noindent I claim that

\[ <u_j^{(k)}(x^{(k)}, B^{(k)}, E^{(k)}, F^{(k)}, G_h^{(k)})u_j^{(k) *}
>_{(h,j)  \in H_k \times J_k}\]

\noindent is an $10 \epsilon \sqrt{6n+3l+2}$-cover for
$\Gamma_R(s,B,E,F,G;m,k,\gamma)$ with respect to $|\cdot|_2.$

\vspace{.1in}

Towards this end suppose $(x,T,P,Q,Z) \in \Gamma_R(s,B,E,F,G;m,k,\gamma)$
for $k$ sufficiently large so that $(x^{(k)}, B^{(k)},E^{(k)},F^{(k)})$
exists as arranged in the preceding paragraph.  By the first paragraph
there exists a $u \in U_k$ such that

\[ |u(x^{(k)}, B^{(k)},E^{(k)}, F^{(k)}) u^* - (x,T,P,Q)|_2 < 
\frac{\epsilon}{(R+1)^2}.\]

\noindent Set $P= \{ p_1,\ldots,p_n\}, Q = \{q_1,\ldots,q_n\}, Z=
\{z_1,\ldots, z_n\}$ and observe that by if $m$ and $\gamma$ are chosen
appropriately, we easily have for each $i$ the inequalities $|z_i|_2 < 2$ 
and $|z_i - p_iz_iq_i|_2 < \epsilon$ and in the situation where $e_i =f_i$ 
the additional condition that $|z_i - (z_i +z_i^*)/2|_2 < \epsilon.$  If $e_i$ 
and $f_i$ are orthogonal set $a_i = z_i$ and otherwise set $a_i = (z_i 
+z_i^*)/2.$  We have:

\[ |z_i - ue_i^{(k)}u^*a_iuf_i^{(k)}u^*|_2 \leq |p_i - ue_i^{(k)}u^*|_2 
\cdot |z_iq_i|_{\infty} + |ue^{(k)}u^*z_i|_{\infty} \cdot |q_i - 
uf_i^{(k)} u^*|_2 + 2 \epsilon < 4\epsilon. \]

\noindent There is an $h \in H_k$ for which $|G^{(k)}_h - (e_1^{(k)} u^*
a_1 u f_i^{(k)}, \ldots, e_n^{(k)} u^* a_n u f_n^{(k)})|_2 < 3 \epsilon
\sqrt{n}.$  By the above $|Z - uG^{(k)}_hu^*|_2 < 7 \epsilon
\sqrt{n}.$ Consequently $|(x,T,P,Q,Z) -
u(x^{(k)},B^{(k)},E^{(k)},F^{(k)},G^{(k)}_h)u^*|_2 < 8 \epsilon \sqrt{n}.$ 
There exists a $j \in J_k$ for which $|u - u_j|_{\infty} < \epsilon,$ so 
that using the fact that $|yz|_2 \leq |y|_{\infty} \cdot |z|_2$ we have

\[|(x,T,P,Q,Z) - u_j (x^{(k)},B^{(k)},E^{(k)},F^{(k)},G^{(k)}_h )u^*_j|_2 < 
10 \epsilon \sqrt{6n+3l+2}\] 

\noindent as promised.

\vspace{.1in}

We've just produced for any $R >0$ and $0 < \epsilon <\epsilon_0$
corresponding $m$ and $\gamma$ such that 

\[\mathbb K_{20 \epsilon
\sqrt{6n+3l+2},R }(s,B,E,F,G;m,\gamma)) < r \cdot (|\log \epsilon| + 
(C+3)).\]

\vspace{.08in}

\noindent Thus, $\mathbb P_{40 \epsilon \sqrt{6n+3l+2}}(s,B,E,F,G) \leq
\mathbb K_{20 \epsilon \sqrt{6n+3l+2}} (s,B,E,F,G) \leq r \cdot (|\log
\epsilon| + (C+3))$ from which it follows that $\mathbb P^r(s,B,E,F,G) <
\infty.$ \end{proof}

By [3] we have:

\begin{corollary} $\mathbb H(s,B,E,F,G) \leq \delta_0(s,B,E,F,G) \leq r.$  
\end{corollary}

\section{Lower Bound}

In this section we show the other inequality, i.e., that the microstate
spaces are rich enough so that their free Hausdorff dimension is $r$ and
we show under additional conditions that their free Hausdorff $r$ entropy
is finite.  The lower bound is more involved.

\vspace{.1in}

Throughout $\kappa(k,r)$ denotes the volume of the ball of radius $r$ in
$\mathbb R^k,$ i.e.,

\[ \kappa(k,r) = \frac{ (r\sqrt{\pi})^k}{\Gamma(\frac{k}{2}+1)}.\]

\vspace{.1in}

$E$ and $F$ generate a commutative, finite dimensional von
Neumann algebra $A$ and thus there exist mutually orthogonal projections
$p_1,\ldots,p_d \in M$ whose span is exactly $A.$ Choose a single
contraction $z=z^*$ which generates $A.$  Observe that $\vee_{j=1}^d p_j = 
\vee_{i=1}^n (e_i\vee f_i).$  Recall the constants $D$ and 
$L$ corresponding to $E\bigcup F, z, R=2,$ and
$p$ as in Section 5 of [3].

\vspace{.1in}

In order to show that generators for finite dimensional algebras
always have nondegenerate free Hausdorff entropy [3] (Lemma 5.3) resorted
to strengthening the work in [2] dealing with representations of
finite dimensional algebras.  This "strengthening" involves replacing
inequalities of the type $\alpha \geq \beta - \epsilon$ for all
$\epsilon>0$ with inequalites of the type $\alpha > \beta.$ We have to do
something similar here under certain conditions on the $e_i$ and $f_i,$ 
but the situation is slightly different.  We will
need to find such bounds not for the dimensions of 
the finite dimensional algebra, but for the dimension of the 
$G$ component (dimensions in the asymptotic matricial microstate sense).
Below we also deal with our general situation where the join of the $e_i$ 
and $f_i$ is not strictly less than the identity.

\begin{lemma} If $1 > \varepsilon > 0$ then there exists an $N \in \mathbb 
N$ such that for each $k > N,$ there is a corresponding 
$*$-isomorphism $\pi_k: A \rightarrow M_k(\mathbb C)$ satisfying:

\begin{itemize} 
\item $\|tr \circ \pi_k - \varphi\| < \varepsilon.$
\item If $H_k$ is the set of unitaries of $\pi_k(A)^{\prime}$ then $H_k$ 
is a tractable subgroup.
\item $\sum_{i=1}^n m_i tr_k(\pi_k(e_i)) tr_k(\pi_k(f_i)) > \sum_{i=1}^n 
m_i \varphi(e_i)\varphi(f_i) -\varepsilon.$
\end{itemize}

\noindent We can arrange it so that if $I \in A,$ then $\pi$ is unital.  
If we assume that $\vee_{j=1}^d p_j < I,$ then the 
third item above can be replaced with the condition $\sum_{i=1}^n m_i 
tr_k(\pi_k(e_i)) tr_k(\pi_k(f_i)) > \sum_{i=1}^n
m_i \varphi(e_i)\varphi(f_i).$

\end{lemma}

\begin{proof} Set $r_j = \varphi(p_j)$ and $t= r_1 + \cdots + r_d.$  
Because $A$ is equal to the span of $\{p_1,\ldots,p_d\}$ there exists 
an $f:\mathbb R^d \rightarrow \mathbb R$ of the form 
$f(t_1,\ldots,t_d) = \sum_{v=1}^L t_{i_v} t_{j_v}$ where $1 
\leq i_1,\ldots,i_L, j_1,\ldots, j_L \leq d,$ and such that for any state 
$\psi$ on $A$ 

\[f(\psi(p_1),\ldots,\psi(p_d)) = \sum_{i=1}^n m_i \psi(e_i) \psi(f_i).
\]

\noindent Suppose $\gamma >0.$ There exist $0 \leq y_1,\ldots,y_d$ such
that $y_1 + \ldots + y_d \leq 1$ and $f(y_1,\ldots,y_d) >
f(r_1,\ldots,r_d) - \varepsilon = \sum_{i=1}^n m_i tr_k(e_i) tr_k(f_i) -
\varepsilon.$ For each $k$ choose $m_1^{(k)},\ldots,m_d^{(k)} \in \mathbb
N \bigcup \{0\}$ such that $|m^{(k)}_i/k - y_i| \leq dk^{-1}$ and
$m^{(k)}_1 + \ldots + m^{(k)}_d \leq k$ with equality in the second
statement if $I \in A.$ Find a representation $\pi_k : A \rightarrow
M_k(\mathbb C)$ such that $tr_k(\pi_k(p_j)) = m_j^{(k)}/k.$ It can easily
be arranged so that $H_k,$ the group of unitaries of $\pi_k(A)^{\prime},$
is tractable.  It is also clear that if $\gamma >0$ is chosen
appropriately from the get-go, then for sufficiently large $k$ $\|tr_k
\circ \pi_k - \varphi \| < \varepsilon.$ Hence we have the first two
conditions and the one in the situation that $I \in A.$  For the last one 
as $k \rightarrow \infty$

\[ \sum_{i=1}^n m_i tr_k(\pi_k(e_i))tr_k(\pi_k(f_i)) = f(m_1^{(k)}, 
\ldots, m_d^{(k)}) \rightarrow f(y_1,\ldots, y_d) >
\sum_{i=1}^n m_i \varphi (e_i) \varphi (f_i) - \varepsilon.\]     

\noindent Finally, if we assume that $\vee_{j=1}^d p_j < I,$ 
then it follows that $0<t<1.$  Thus, we can repeat the same argument 
except choosing the $y_i$ to satisfy $y_1 + \cdots + y_d  \leq 1$ and 
$f(y_1,\ldots,y_d) > f(r_1,\ldots,r_d).$  We again get the first two 
conditions and the line above, except that we don't have the 
$-\varepsilon.$ 

\end{proof}

\begin{lemma} $\mathbb H(s,B,E,F,G) =r$ and if $\hspace{.06in} 
\vee_{j=1}^d p_j <I,$ then $- \infty < \mathbb H^r(s,B,E,F,G) < \infty.$
\end{lemma}

\begin{proof} By Lemma 2.1 and Corollary 2.2 we have $\mathbb H(s,B,E,F,G)
\leq r$ and $\mathbb H^r(s,B,E,F,G) < \infty$ so it suffices to prove the
other inequality for both statements.  Suppose $m \in \mathbb N$ and
$t, \gamma >0$ are given.  There exist $m_1 \in \mathbb N$ and $\gamma_1 
>0$ such that if $(\eta, Y_0, Y_1, Y_2) \in 
\Gamma_2(s,B,E,F;m_1,k,\gamma_1)$
and $Y_3 \in \Gamma_2(s_1, \ldots, s_n;m_1,k,\gamma_1)$ are
$(m_1,\gamma_1)$-free, then

\[ (\eta, Y_0, Y_1, Y_2, Y_1 Y_3 Y_2) \in
\Gamma_8(s,B,E,F,G;m,k,\gamma)\] 

\noindent where $Y_1 Y_3 Y_2$ is the $n$-tuple obtained from multiplying
the entries of the three $n$-tuples in the natural way.  Also by [2] and
the existence of finite dimensional approximants for any tuple which
generates a hyperfinite algebra, there exist $m_2 \in \mathbb N$ and
$\gamma_2 >0$ such that for $k$ sufficiently large if $(\eta, Y_1, Y_2)  
\in \Gamma_2(s,E,F;m_2,k,\gamma_2),$ then there exists an $l$-tuple $Y_0$
such that $(\eta, Y_0,Y_1, Y_2) \in \Gamma_2(s,B,E,F;m_1,k,\gamma_1).$

\vspace{.1in}

Define $\delta = 0$ if $ \vee_{j=1}^d p_j < I$ and $\delta = t$ otherwise.  
Replacing Lemma 3.6 of [2] with Lemma 3.1 above in the proof of Lemma 5.2
of [2] (which works in the nonunital case) produces $1 > \lambda, \zeta, c
>0$ independent of the $m$ and $\gamma$ such that for sufficiently large
$k$ there exists a $*$-homomorphism $\pi_k : A \rightarrow M_{k}(\mathbb
C)$ which is unital if $I \in A$ and:

\begin{itemize}

\item $\|tr_k \circ \pi_k - \varphi \|
< \frac {\gamma_2}{2^{m_2}}.$

\vspace{.05in}

\item Denote by $H_k$ the set of unitaries of $\pi_k(A)^{\prime}.$
Define $\mathcal H_k \subset iM^{sa}_k(\mathbb C)$ to be the Lie
subalgebra of $H_k$ and $\mathcal X_k$ to be the orthogonal
complement of $\mathcal H_k$ with respect to the Hilbert-Schmidt inner
product.  For every $s>0$ write $\mathcal X^s_k$ for the ball in $\mathcal
X_k$ of operator norm less than or equal to $s.$  Below all
volume quantities are obtained from Lebesgue measure when the spaces are
given the inner product $Re Tr.$

\[ \frac{\text{vol}(\mathcal X^1_k)}{C(\dim \mathcal X_k, \sqrt{k})} > 
(\zeta)^{\dim \mathcal X_k}.\]

\vspace{.05in}

\item For any $x,y \in \mathcal X^c_k$
\[d_2(q(e^x), q(e^y)) \geq \lambda |x-y|_2.\]

\vspace{.05in}

\item $\sum_{i=1}^n m_i \cdot tr_k(\pi_k(e_i)) \cdot
tr_k(\pi_k(f_i)) > \sum_{i=1}^n m_i \varphi(e_i) \varphi(f_i) - \delta.$

\end{itemize}

\vspace{.1in}

	For $k$ sufficiently large and $1 \leq j \leq d$ set $p_j^{(k)} =
\pi_k(p_j),$ $p_{d+1}^{(k)} = I - (p^{(k)}_1+ \cdots + p_d^{(k)})$
($p_{d+1}^{(k)} = 0$ if $I \in A$), and for $1 \leq i \leq n$ set
$e^{(k)}_i = \pi_k(e_i)$ and $f^{(k)}_i = \pi_k(f_i).$ Define $E^{(k)} =
\{e_1^{(k)}, \ldots, e_n^{(k)}\}$ and $F^{(k)} = \{f_1^{(k)},\ldots,
f_n^{(k)}\}.$ For each $1 \leq j \leq d, sp_j \in (p_jMp_j,
\varphi(p_j)^{-1} \cdot \varphi)$ has finite $\chi^{sa}$ entropy and in 
this same way so does $sp_{d+1}$ if $I \notin A.$  Define

\[ X_k = \bigoplus_{j=1}^{d+1} \Gamma_2(sp_j;m_2,tr_k (\pi_k(p_j))k,
\gamma_2) \subset \bigoplus_{j=1}^{d+1} p_j^{(k)}M^{sa}_k(\mathbb 
C)p^{(k)}_j \subset M^{sa}_k(\mathbb C).\]

\noindent Write $\lambda_k$ for the measure on $M^{sa}_k(\mathbb C)$ given
by Lebesgue measure on $\bigoplus_{i=1}^{d+1} p_i M^{sa}_k(\mathbb C)p_i$ 
with respect to the $\| \cdot \|_2$ norm restricted to this direct sum.  
Define $Y_k = \Gamma_2(s_1,\ldots,s_n;m_1,k,\gamma_1)$ and write $\mu_k$
for Lebesgue measure on $(M^{sa}_k(\mathbb C))^n$ with respect to the
$\|\cdot\|_2$ norm on $(M^{sa}_k(\mathbb C))^n.$ Consider the probability
measure $\sigma_k$ on $((M^{sa}_k(\mathbb C))_{2})^{n+1}$ obtained by
restricting $\lambda_k \times \mu_k$ to $X_k \times Y_k$ and normalizing
appropriately.  Lemma 2.14 of [9] provides an $N \in \mathbb N$ such that
if $k \geq N$ and $\sigma$ is any Radon probability measure on
$((M^{sa}_k(\mathbb C))_2)^{3n+l+1}$ invariant under the $U_k$-action

\[ (\xi_1,\ldots, \xi_{3n+l+1}) \mapsto (\xi_1,\ldots, 
\xi_{2n+l+1}, u\xi_{2n+l+2}u^*, \ldots, u\xi_{3n+l+1}u^*)
\] 

\vspace{.1in}

\noindent then $\sigma(\omega_k) > \frac{1}{2}$ where $\omega_k$ is the
subset of $((M_k^{sa}(\mathbb C))_2)^{3n+l+1}$ consisting of those tuples
such that the first $2n+l+1$ entries are $(m_1, \gamma_1)$-free from the
last $n.$ By the first paragraph if $k$ is large enough for each $\xi \in
X_k,$ then there exists an $l$-tuple $T$ such that $(\xi
,T,E^{(k)},F^{(k)}) \in \Gamma_2(s,B,E,F;m_1,k,\gamma_1).$ Define
$\delta_{\xi}$ to be the atomic probability measure on $((M_k^{sa}(\mathbb
C))_2)^{2n+l+1}$ supported at $(\xi,T,E^{(k)},F^{(k)}).$ Writing
$\overline{\mu_k}$ for the normalization of $\mu_k$ we have that
$\delta_{\xi} \times \overline{\mu_k}$ is a Radon probability measure on
$((M_k^{sa}(\mathbb C))_2)^{3n+l+1}$ invariant under the $U_k$-action
described above so that $(\delta_{\xi}\times \overline{\mu_k})(\omega_k) >
\frac{1}{2}.$ Define $\Theta_k$ to be the set of all $(n+1)$-tuples
$(\xi_1,\ldots,\xi_{n+1})$ for which:

\begin{itemize}

\item $\xi_{1} \in X_k$ and $(\xi_{2}, \ldots, \xi_{n+1} )  
\in Y_k.$

\item  There is a $T$ satisfying the two conditions that
$(T, E^{(k)}, F^{(k)},\xi_1) \in 
\Gamma_2(B,E,F,s;m_1,k,\gamma_1)$ 
and  $(T,E^{(k)},F^{(k)},\xi_1,\ldots,\xi_{n+1}) \in
\omega_k.$ \end{itemize}

\noindent $\Theta_k \subset X_k \times Y_k$ is an open (and thus 
measurable) set.  The fact that $(\delta_{\xi} 
\times\overline{\mu_k})(\omega_k) > 1/2$ for every $\xi \in X_k$ in conjunction 
with Fubini's Theorem tells us that $\sigma_k(\Theta_k) > \frac{1}{2}.$

\vspace{.15in}
   
On $V_k = (\bigoplus_{j=1}^{d+1} p^{(k)}_jM^{sa}_k(\mathbb C)p^{(k)}_j) 
\bigoplus (\bigoplus_{i=1}^n M_k^{sa}(\mathbb C)) \subset 
(M^{sa}_k(\mathbb
C))^{n+1}$ consider the real orthogonal projection $Q_k$ on
$(M^{sa}_k(\mathbb C))^{n+1}$ defined by

\[ Q_k(\xi, \eta_1,\ldots,\eta_n) = (\xi, r_1(e_1^{(k)} \eta_1
f_1^{(k)} + f_1^{(k)} \eta_1 e_1^{(k)}), \ldots, r_n(e_n^{(k)}
\eta_n f_n^{(k)} + f_n^{(k)} \eta_n e_n^{(k)}))  \]
  
\vspace{.1in}

\noindent where $r_i = 1/(3 - m_i).$ Denote $a_k = \dim V_k - \dim
Q(V_k). $ If $m_k$ is Lebesgue measure on $Q_k(V_k)$ obtained with respect
to the $\| \cdot \|_2$ norm of $(M^{sa}_k(\mathbb C))^{n+1}$ restricted to
$Q_k(V_k),$ then

\[ \frac{1}{2} \cdot (\lambda_k \times \mu_k)(X_k \times
Y_k)  < (\lambda_k \times \mu_k)(\Theta_k) \leq  \kappa (a_k, 
\sqrt{2nk} )  \cdot m_k(Q_k(\Theta_k)). \]

\vspace{.1in}

\noindent Thus, $m_k(Q_k(\Theta_k)) \geq \frac{1}{2} \cdot \lambda_k(X_k) 
\lambda_k(Y_k) \cdot \kappa(a_k, \sqrt {2nk} )^{-1}.$  Define the linear 
map $P_k$ from $(M^{sa}_k(\mathbb C))^{n+1}$ into $\bigoplus_{i=1}^{n+1} 
M_k(\mathbb C)$ by 

\[ P_k(\xi, \eta_1,\ldots, \eta_n) = (\xi, e_1^{(k)}\eta_1 f_1^{(k)}, 
\ldots, e_n^{(k)} \eta_n f_n^{(k)}).\]

\noindent $P_k$ is the the composition of $Q_k$  
with a bi-Lipschitz map bounded from below by $1/2$ and above by $1,$ this 
bi-Lipschitz map defined by  

\[Q_k(\xi, \eta_1, \ldots, \eta_n) \mapsto (\xi, e_1^{(k)} \eta_1 
f_1^{(k)}, \ldots, e_n^{(k)} \eta_n f_n^{(k)}).\]

\noindent Consequently, if we endow the range of $P_k$ with the inherited 
$\| \cdot \|_2$-norm the Lebesgue measure of $\Theta_k$ (with respect to 
this identification) is no less than 

\[ \frac{1}{2^{\dim P_k}} \cdot  \lambda_k(X_k)
\lambda_k(Y_k) \cdot \kappa (a_k, \sqrt{2nk})^{-1}.\] 

\vspace{.15in}

Define $\Omega_k$ to be the set of all elements of the form 
$u(T,E^{(k)},F^{(k)}, Y)u^*$ where $u = e^x, x \in \mathcal 
X^r_k,$ $Y \in P_k(\Theta_k),$ and $(T,E^{(k)},F^{(k)},Y) \in 
\Gamma_2(B,E,F,s,G;m,k,\gamma).$  Clearly $\Omega_k \subset 
\Gamma_2(B,E,F,s,G;m,k,\gamma).$  Consider the map $\Phi : \Omega_k 
\rightarrow \mathcal X^c_k \times 
P_k(\Theta_k)$ defined by 

\[ \Phi(u(T, E^{(k)},F^{(k)}, Y)u^*) = (x, Y)
\]

\noindent where $u = e^x$ for some $x \in \mathcal X^c_k$ and $Y \in 
P_k(H_k).$  This map is well-defined for suppose $u=e^x, v  = 
e^{x^{\prime}}$ for some $x, x^{\prime} \in \mathcal X^c_k,$ 
and $u(T, E^{(k)},F^{(k)}, Y)u^* = 
v(T,^{\prime},E^{(k)},F^{(k)}, Y^{\prime})v^*$ where 

\[(T, E^{(k)},F^{(k)}, 
Y),(T,^{\prime},E^{(k)},F^{(k)}, Y^{\prime}) \in 
\Gamma_2(B,E,F,s,G;m,k,\gamma).\]

\noindent $v^*u(E^{(k)},F^{(k)})u^*v = (E^{(k)},F^{(k)}).$  By definition 
$v^*u \in H_k$ so that 

\[ 0 = d_2(q(e^x),q(e^{x^{\prime}})) \geq 
\lambda |x - x^{\prime}|_2.\]  

\noindent $x = x^{\prime} \Rightarrow u= v \Rightarrow y = y^{\prime}.$
$(x,Y) = (x^{\prime},Y^{\prime})$ and thus, $\Phi$ is well-defined.  
$\Phi$ is also Lipschitz for suppose that $u = e^x, v= e^{x^{\prime}}$ for
$x,x^{\prime} \in \mathcal X^c_k,$  $Y, Y^{\prime} \in L(\Theta_k),$ and
there exist $T, T^{\prime}$ for which $(T,E^{(k)},F^{(k)},Y),
(T^{\prime},E^{(k)},F^{(k)},Y^{\prime}) \in
\Gamma_2(B,E,F,s,G;m,k,\gamma).$

\[
|\Phi(u(T,P^{(k)},Q^{(k)},Y)u^*) - 
\Phi(v(T^{\prime},P^{(k)},Q^{(k)},Y^{\prime})v^*)|_2 = 
|(x,Y) - (x^{\prime},Y^{\prime})|_2 \leq  |x - x^{\prime}|_2 + |Y - 
Y^{\prime}|_2.
\]
\noindent The analysis of Lemma 5.4 in [3] shows that 
there exist constants $D$ and $L$ dependent only on $A, e_1, \ldots, e_n, 
f_1,\ldots,f_n,z,$ and $p$ such that

\[ 
|x - x^{\prime}|_2 \leq \frac{DL}{\lambda} \cdot |u(P^{(k)}, Q^{(k)})u^* - 
v(P^{(k)},Q^{(k)})v^*|_2
\]

\noindent and

\begin{eqnarray*} |Y - Y^{\prime}|_2 & \leq & |uYu^* - vY^{\prime}v^*|_2 + 
|vY^{\prime}v^* - uY^{\prime}v^*|_2 + |uY^{\prime}v^* - uY^{\prime}u^*|_2 
\\ & \leq & |uYu^* -vY^{\prime}v^*|_2 + 4|e^x-e^{x^{\prime}}|_2 \\
& \leq & |uYu^* -vY^{\prime}v^*|_2 + 4|x- x^{\prime}|_2 \\
& \leq & |uYu^* -vY^{\prime}v^*|_2 + \frac{4DL\sqrt{n+1}}{\lambda} \cdot  
|u(P^{(k)},Q^{(k)})u^*-v(P^{(k)},Q^{(k)})v^*|_2.\\   
\end{eqnarray*}

\noindent From this it follows that $\| \Phi \|_{Lip} \leq C$ where $C =
6(DL\sqrt{n+1}+1)\lambda^{-1}.$ Finally, the range of $\Phi$ is exactly
$\mathcal X^r_k \times P_k(\Theta_k).$ This follows from the way in which 
we defined $\Theta_k$.  For given $x \in \mathcal X^c_k$ and $Y \in
P_k(\Theta_k)$ there exists some $(z_1,\ldots,z_{n+1})=Z \in \Theta_k$ for
which $P_k(Z) = Y.$ By definition of $\Theta_k$ we have
$(z_2,\ldots,z_{n+1})\in \Gamma_2(s_1,\ldots,s_n;m_1,k,\gamma_1),$ the
existence of a $T$ for which $(T, E^{(k)},F^{(k)}, z_1) \in
\Gamma_2(B,E,F,s;m_1,k,\gamma_1),$ and $(T,E^{(k)},F^{(k)},Z) \in
\omega_k.$ By the first paragraph \[(T,E^{(k)},F^{(k)},P_k(Z)) =
(T,E^{(k)},F^{(k)},Y) \in \Gamma_2(B,E,F,s,G;m,k,\gamma).\]

\noindent Consequently, $\Phi(e^x(T, E^{(k)},F^{(k)}, Y)e^{-x}) = (x, 
Y)$ as desired.    

\vspace{.1in}

The preceding paragraph shows that

\[ C^{rk^2} \cdot H^{rk^2}_{\epsilon}(\Gamma_2(B,E,F,s,G;m,k,\gamma)) \geq 
H^{rk^2}_{C \epsilon}(\mathcal X^c_k \times P_k(\Theta_k)) 
\]

\noindent So we just to need to approximate the right hand side above and
we will do so using by comparing volumes.  Suppose $0 < \epsilon <
C^{-1}.$ Suppose $<\theta_j>_{j \in J}$ is a countable $C
\epsilon$-cover of $\mathcal X^c_k \times P_k(\Theta_k).$ Without loss of
generality we may assume the $\theta_j$ are closed.  Regarding $\mathcal
X^c_k \subset \mathcal X_k,$ the third condition imposed on the $\pi_k$
says 

\[ \frac{\text{vol}(\mathcal X_k^1)}{\kappa(\dim \mathcal X_k, \sqrt{k})}
> (\zeta)^{\dim \mathcal X_k}.\]

\noindent It follows that $\mathcal X_k^c \times P_k(\Theta_k)$ is subset
of $\mathcal X_k \times P_k(M^{sa}_k(\mathbb C))^{n+1}$ with Lebesgue
volume (again computed when the ambient space is endowed with 
$\| \cdot \|_2$) no less than

\[ (c \zeta)^{\dim \mathcal X_k} 
\cdot 2^{- \dim P_k} \cdot \lambda_k(X_k) \mu_k(Y_k) \cdot \frac{ 
\kappa(\dim \mathcal X_k, \sqrt{k}) }{\kappa(a_k, 
\sqrt{2nk})}.\] 

\noindent Also observe that if $b_k$ denotes the dimension of $\mathcal 
X_k \times P_k(M^{sa}_k(\mathbb C))^{n+1},$ then 

\begin{eqnarray*} b_k & = & k^2 \left [ 1 - \sum_{j=1}^{d+1} 
tr_k(p^{(k)}_j)^2 
+\sum_{i=1}^{d+1} 
tr_k(p_j^{(k)})^2 + \sum_{i=1}^n m_i tr_k(e_i^{(k)}) tr_k(f_i^{(k)}) 
\right] \\
& = & k^2 \left [1 +  \sum_{i=1}^n m_i tr_k(e_i^{(k)}) 
tr_k(f_i^{(k)})\right] \\ & > & k^2 (1 - \sum_{i=1}^n m_i \alpha_i 
\beta_i - \delta) \\ & = & (r- \delta)k^2.\\
\end{eqnarray*}

\noindent Thus, using the preceding volume estimates with the lower bound 
on $b_k$ we have 

\[ \sum_{j \in J} |\theta_j|^{(r- \delta)k^2} \geq \sum_{j \in J}
|\theta_j|^{b_k \cdot k^2}  \geq (c 
\zeta)^{\dim \mathcal X_k}
\cdot 2^{- \dim P_k} \cdot \lambda_k(X_k) \mu_k(Y_k) \cdot \frac{     
\kappa(\dim \mathcal X_k, \sqrt{k}) }{\kappa(a_k,
\sqrt{2nk})  \kappa(b_k, \sqrt{k})  } 
\]

\noindent It now remains to compute the asymptotics of the right hands
side.  Suppose $0 <\epsilon < C^{-1}$ and $m$ and $\gamma$ are given.  
Regularity of a single semicircular, regularity of a free family of
semicirculars, $\chi^{sa}(sp_j) > -\infty,$ and Stirling's Formula imply
that $\limsup_{k\rightarrow \infty} k^{-2} \cdot \log
H^{(r-\delta)k^2}_{\epsilon}(\Gamma_2(B,E,F,s,G;m,k,\gamma))$ dominates

\begin{eqnarray*} & & \limsup_{k \rightarrow \infty} k^{-2} \cdot \log 
\left( (c \zeta)^{\dim \mathcal 
X_k} \cdot 2^{- \dim P_k} \cdot \lambda(X_k) \cdot \mu_k(Y_k) \cdot 
\frac{\kappa(\dim \mathcal X_k, \sqrt{k})}{\kappa(a_k,\sqrt{2nk}) 
\kappa(b_k,\sqrt{k})} \cdot C^{-rk^2} \right) \\ & \geq &  
\limsup_{k \rightarrow \infty} k^{-2} \log \left(\Pi_{j=1}^d 
\kappa(tr_k(\pi_k(p_j))^2 k^2, \sqrt{k}) \cdot (\kappa(k^2, \sqrt{k})^n 
\cdot (\kappa(a_k, \sqrt{k})^{-1} \cdot (\kappa(b_k, \sqrt{k})) ^{-1} 
\right) \\ & & + (n+1) \log \left(\frac{c\zeta}{2C} \right) \\ & > & - 
\infty. \end{eqnarray*}

\noindent The estimate hold for arbitrary $m$ and $\gamma$ and $\epsilon$
sufficiently small whence $\mathbb H^{r-\delta}(B,E,F,s,G) = \mathbb H^{r-
\delta}(s,B,E,F,G)  > -\infty.$ Now in the general case this holds for all
$\delta = t>0$ and thus we have that $\mathbb H(s,B,E,F,G) \geq r.$ If
$\vee_{j=1}^d p_j <I,$ then $\delta =0$ and we arrive at
$\mathbb H^r(s,B,E,F,G) > -\infty.$ \end{proof}

\vspace{.1in}

\noindent{\it Acknowledgements.} I thank Dan Voiculescu, my advisor, who
suggested this problem and an alternative interpretation of $\mathbb H^r.$
Part of this research was conducted at Paris VII and I thank Georges
Skandalis and the operator algebra team for their hospitality.

\end{document}